\input amstex\documentstyle{amsppt}  
\pagewidth{12.5cm}\pageheight{19cm}\magnification\magstep1
\topmatter
\title Traces on Iwahori-Hecke algebras and counting rational points\endtitle
\author G. Lusztig\endauthor
\address{Department of Mathematics, M.I.T., Cambridge, MA 02139}\endaddress
\thanks{Supported in part by NSF grant DMS-1855773.}\endthanks
\endtopmatter   
\document
\define\tfB{\ti{\fB}}
\define\tch{\ti{\ch}}

\define\Irr{\text{\rm Irr}}

\define\frl{\forall}

\define\qua{\quad}

\define\lb{\linebreak}

\define\part{\partial}

\define\m{\mapsto}
\define\do{\dots}

\define\bsl{\backslash}

\define\sub{\subset}    

\define\T{\times}
\define\ti{\tilde}
\define\nl{\newline}
\redefine\i{^{-1}}

\define\un{\underline}
\define\ov{\overline}
\define\ot{\otimes}
\define\bbq{\bar{\QQ}_l}

\define\Hom{\text{\rm Hom}}
\define\End{\text{\rm End}}

\define\Ind{\text{\rm Ind}}

\define\sgn{\text{\rm sgn}}
\define\tr{\text{\rm tr}}

\redefine\d{\delta}

\define\p{\pi}

\define\r{\rho}

\redefine\t{\tau}

\redefine\l{\lambda}
\define\z{\zeta}
\define\x{\xi}

\redefine\G{\Gamma}
\redefine\D{\Delta}

\define\kk{\bold k}

\define\qq{\bold q}

\define\FF{\bold F}

\define\NN{\bold N}

\define\QQ{\bold Q}

\define\ZZ{\bold Z}

\define\cb{\Cal B}

\define\ce{\Cal E}
\define\cf{\Cal F}

\define\ch{\Cal H}

\define\cl{\Cal L}

\define\co{\Cal O}

\define\ct{\Cal T}

\define\fB{\frak B}

\define\sha{\sharp}

\head Introduction\endhead
\subhead 0.1\endsubhead
Let $W$ be a Weyl group with length function $w\m|w|$ and let $\ch$ be the Iwahori-Hecke algebra
over $\QQ(\qq)$ ($\qq$ is an indeterminate) attached to $W$. Recall that $\ch$ is the
$\QQ(\qq)$-vector space with basis $\{T_w;w\in W\}$ with multiplication given by
$T_wT_{w'}=T_{ww'}$ if $w,w'\in W$ satisfy $|ww'|=|w|+|w'|$ and $(T_t-\qq)(T_t+1)=0$ if $t\in W,|t|=1$.
For $w\in W$ let $\t(w,\qq)\in\ZZ[\qq]$ be the trace of the linear map $\ch@>>>\ch$ given by $h\m T_wh$.
The ``trace polynomials''
$\t(w,\qq)$ appear in relation with counting $\FF_q$-rational points in certain algebraic
varieties; three apparitions are in \cite{L78},\cite{L21},\cite{L85}, see 0.2(a), 0.3(a),
0.4(a); the fourth one is new, see 0.5.
(We denote by $\FF_q$ the subfield with $q$ elements of an algebraic closure $\kk$ of $\FF_p$; $q$ is a
power of a prime number $p$.)

\subhead 0.2\endsubhead
Let $G$ be a connected, reductive group over $\kk$ whose Weyl group
is $W$. We fix a split $\FF_p$-rational structure on $G$ with Frobenius map $F:G@>>>G$.
Let $\cb$ the variety of Borel subgroups of $G$. Then $W$ can be viewed as the set of orbits on $\cb\T\cb$
for the diagonal conjugation $G$-action. For any $w\in W$ let $\co_w\sub\cb\T\cb$ be the corresponding
$G$-orbit. Let $F:\cb@>>>\cb$ be the map induced by $F:G@>>>G$.
For any $B\in\cb$ let $U_B$ be the unipotent radical of $B$.
We fix $B^+\in B$ such that $F(B^+)=B^+$.
We set $\nu=\dim\cb$. Let $r$ be the rank of $G$.

For $w\in W,s\ge1$ let
$X_{w,s}=\{B\in\cb;(B,F^s(B))\in\co_w\}$ (see \cite{DL}). This is a subvariety of $\cb$ on which
$G^{F^s}$ acts by conjugation; hence there is an induced action of $G^{F^s}$ on the $\bbq$-cohomology space $H^i_c(X_{w,s})$ so that the virtual representation
$R^1_{w,s}=\sum_i(-1)^iH^i_c(X_{w,s})$ of $G^{F^s}$ is defined.  
(For a set $Z$ and a map $F:Z@>>>Z$ we write $Z^F=\{z\in Z;F(z)=z\}$.)

Let $s'\ge1$. The following result is obtained from of \cite{L78, 3.8} with $w'=1$.

(a) {\it The number of $\FF_{p^{ss'}}$-rational points of $(B^+)^{F^s}\bsl X_{w,s}$
(where $(B^+)^{F^s}$ acts by conjugation) is equal to $\t(w,p^{ss'})$.
Moreover, we have $H^i_c((B^+)^{F^s}\bsl X_{w,s})_j=0$ for $j$ odd and
$\sum_{i,j}(-1)^i\dim H^i_c((B^+)^{F^s}\bsl X_{w,s})_j\qq^{j/2}=\t(w,\qq)$.}
\nl
(Here $()_j$ is the part of pure weight $j$.)

\subhead 0.3\endsubhead
For $w\in W,g\in G$ let $Y_{g,w}=\{B\in\cb;(B,gBg\i))\in\co_w\}$.
Let $s\ge1$ and let $g\in(B^+)^{F^s}$ be regular semisimple in $G$.
Then $Y_{g,w}$ is defined over $\FF_{p^s}$.
The following result can be found in \cite{L21, 0.7(a),0.8(a)}, see also \cite{L79, 1.2}.

(a) {\it The number of $\FF_{p^s}$-rational points of $Y_{g,w}$ is equal to $\t(w,p^s)$.
Moreover, we have $H^i_c(Y_{g,w})_j=0$ for $j$ odd and
$\sum_{i,j}(-1)^i\dim H^i_c(Y_{g,w})_j\qq^{j/2}=\t(w,\qq)$.}

\subhead 0.4\endsubhead
Let $w\in W,s\ge1$. Let
$$\tfB_w=\{(g,B,B')\in G\T\cb\T\cb; g\in B,(B',gB'g\i)\in\co_w\}.\tag a$$
This variety is naturally defined over $\FF_{p^s}$. The following statement is
a special case of \cite{L85, 13.5} with $\cl=\bbq,w'=1$:

(b) {\it $(\tfB_w)^{F^s}=\sha(G^{F^s})\t(w,p^s)$.}

\subhead 0.5\endsubhead
Let $\ct$ be ``the maximal torus'' of $G$ identified with $B/U_B$ for any $B\in\cb$.
Let $w\in W$. We define a morphism $\p:\tfB_w@>>>\ct$ by $\p(g,B,B')=$ image of
$g$ in $B/U_B=\ct$. For $\d\in\ct$ we set
$$\tfB_w^\d=\p\i(\d)=\{(g,B,B')\in\tfB_w; g\in B_\d\}$$
where $B_\d$ is the inverse image under $B@>>>B/U_B$ of the element $\d_B\in B/U_B$ defined by $\d$.
The following refinement of 0.4(b) will be proved in 1.2.

(a) {\it If $s\ge1$ and $\d\in\ct^{F_s}$ then}
$\sha((\tfB_w^\d)^{F^s})=p^{\nu s}\sha(\cb^{F^s})\t(w,p^s).$

\subhead 0.6\endsubhead
Specializing $\d=1$, we have
$$\tfB_w^1=\{(g,B,B')\in\tfB_w;g\text{ unipotent}\}.$$
The variety $\tfB_w^1$ is closely related to the main theme of \cite{LY} (see below); it has been
been considered independently by M.-T. Trinh.

Let $\Irr(W)$ be the set of irreducible representations of $W$ (up to isomorphism).
For $w\in W$ we set

$B^+wB^+=\{g\in G;(B^+,gB^+g\i)\in\co_w\}$, 

$(B^+wB^+)^1=\{g\in B^+wB^+;g\text{ unipotent }\}$.
\nl
In \cite{LY, 1.3(b)} for any $w\in W,s\ge1$ and $E\in\Irr(W)$ we have defined
$$\z_{w,p^s,E}=\sum_{g\in((B^+wB^+)^1)^{F^s}}\sum_i(-1)^i\tr(F^s,H^i(\cb_g)_E).\tag a$$
Here $\cb_g=\{B\in\cb;g\in B\}$ and
$H^i(\cb_g)_E=\Hom_W(E,H^i(\cb_g))$ where the Springer
action of $W$ on $H^i(\cb_g)$ is used. We have $\z_{w,p^s,E}\in\ZZ$, see \cite{LY, 1.3(c)}.
Let
$$\widetilde{(B^+wB^+)^1}=\{(g,B)\in(B^+wB^+)^1\T\cb;g\in B\}.$$

From the definition we see that

(b) $\sha(\widetilde{(B^+wB^+)^1})^{F^s})=\sum_{E\in\Irr(W)}\z_{w,p^s;E}\dim(E)$.
\nl
We have

(c) $\sha((\tfB_w^1)^{F^s})=\sha(\cb^{F^s})\sha(\widetilde{(B^+wB^+)^1})^{F^s})
=\sha(\cb^{F^s})\sum_{E\in\Irr(W)}\z_{w,p^s;E}\dim(E)$.
\nl
Indeed, we have a fibration $\tfB_w^1@>>>\cb$, $(g,B,B')\m B'$ whose fibre at $B^+$ is
$\widetilde{(B^+wB^+)^1}$.

Now the numbers $\z_{w,p^s;E}$ are described in \cite{LY} so
that (c) gives a way to compute $\sha((\tfB_w^1)^{F^s})$ other than (and more complicated than)
that given by 0.5(a).

\subhead 0.7\endsubhead
Let $w\in W$. The variety
$$\fB_w=\{(g,B')\in G\T\cb;(B',gB'g\i)\in\co_w\}$$
appeared in the theory of character sheaves. There is a unique map $\bar\p:\fB_w@>>>W\bsl\ct$
such that the diagram
$$\CD\tfB_w@>\p>>\ct\\
    @VVV        @VVV\\
\fB_w@>\bar\p>>W\bsl\ct
\endCD$$
is commutative. 
(Here the left vertical map is $(g,B,B')\m(g,B')$ and the right vertical map is the
obvious one.) For $\d\in\ct$ we denote by $(\d)$ the $W$-orbit of $\d$ and we set
$\fB_w^{(\d)}=\bar\p\i((\d))$. In particular we have
$$\fB_w^{(1)}=\{(g,B')\in G\T\cb;(B',gB'g\i)\in\co_w, g\text{ unipotent}\}.$$
This variety is studied in \S2.

\subhead 0.8\endsubhead
In \S3 we show that, in the case where $w\in W$ is elliptic and of minimal length in its conjugacy
class, the varieties $\tfB_w,\tfB_w^\d$ (with $\d\in\ct$) are smooth, irreducible.
In \S4 we study the space of orbits of the natural $G$-action on $\tfB_w,\fB_w$ with $w$ as above.
In \S5 we prove a symmetry property of the polynomial $\t(w,\qq)$.

\subhead 0.9\endsubhead
I thank Zhiwei Yun, Minh-Tam Trinh and Xuhua He for useful comments.

\head 1. Counting rational points of $\tfB_w^1$\endhead
\subhead 1.1\endsubhead
In this section we fix $s\ge1$. Let $\ch_s$ be the $\bbq$-algebra obtained by replacing $\qq,\QQ[\qq]$ by $p^s,\bbq$
in the definition of
$\ch$ in 0.1. It is known that $\ch_s$ is canonically
isomorphic to the group algebra $\bbq[W]$ provided a square
root
of $p$ is chosen. Hence any $E\in\Irr(W)$ can be regarded as
an irreducible $\ch_s$-module $E_s$. Let $\cf$ be the
vector space of functions $\cb^{F^s}@>>>\bbq$. This is
naturally a $G^{F^s}$-module and we can identify
$\ch_s=\End_{G^{F^s}}\cf$. For $E\in\Irr(W)$ we set
$\cf_{E,s}=\Hom_{\ch_s}(E_s,\cf)$. This is naturally an
irreducible representation of $G^{F^s}$.

\subhead 1.2\endsubhead
We fix $w\in W$ and $\d\in\ct^{F^s}$. Let $N=\sha((\tfB_w^\d)^{F^s})$.

We have $N=\sum_{B\in\cb^{F^s}}\sum_{g\in(B_\d)^{F^s}}\sha(Y_{g,w}^{F^s})$.
By the argument in \cite{L11, 1.3(a)} for any $g\in G^{F^s}$ we have
$$\sha(Y_{g,w}^{F^s})=\sum_{E\in\Irr(W)}\tr(T_w,E_s)\tr(g,\cf_{E,s}).\tag a$$
It follows that 
$$N=\sum_{B\in\cb^{F^s}}\sum_{E\in\Irr(W)}\tr(T_w,E_s)\sum_{g\in(B_\d)^{F^s}}
\tr(g,\cf_{E,s}).$$
For $B\in\cb^{F^s}$ let 
$$\cf_{E,s}^{U_B^{F^s}}=\{\x\in\cf_{E,s};
u\x=\x\qua\frl u\in U_B^{F^s}\}.$$
Note that $\cf_{E,s}^{U_B^{F^s}}$ is stable under the action
of $B^{F^s}$ and is a direct sum of one dimensional
representations of $B^{F^s}$ (on which
$U_B^{F^s}$ acts trivially. If $\l$ is one of these one dimensional representations
we see that $\cf_{E,s}$ appears in
$\Ind_{B^{F^s}}^{G^{F^s}}(\l)$ and also in
$\Ind_{B^{F^s}}^{G^{F^s}}(1)$; but these two induced representations are disjoint unless
$\l=1$. We see that $B^{F^s}$ acts trivially on
$\cf_{E,s}^{U_B^{F^s}}$.
It follows that for any $g_0\in B^{F^s}$ we have
$$\sum_{u\in U_B^{F^s}}\tr(g_0u,\cf_{E,s})=
\sha(U_B^{F^s})\dim(\cf_{E,s}^{U_B^{F^s}}).$$
It is well known that
$$\dim(\r_{E,p^s}^{U_B^{F^s}})=\dim E.$$
It follows that
$$\sum_{u\in U_B^{F^s}}\tr(g_0u,\cf_{E,s})=
\sha(U_B^{F^s})\dim(E).$$
In particular we have
$$\sum_{g\in(B_\d)^{F^s}}\tr(g,\cf_{E,s})
=p^{\nu s}\dim(E).$$
We see that 
$$N=p^{\nu s}\sha(\cb^{F^s})\sum_{E\in\Irr(W)}
\tr(T_w,E_{p^s})\dim(E).$$
Note that $\sum_{E\in\Irr(W)}\tr(T_w,E_s)\dim(E)=\t(w,p^s)$.
(The trace defining $\t(w,p^s)$
can be computed by decomposing the left regular representation
of $\ch_s$ into irreducible submodules.) We see that
$$N=p^{\nu s}\sha(\cb^{F^s})\t(w,p^s).$$

This completes the proof of 0.5(a).

\head 2. Counting rational points of $\fB_w^{(1)}$\endhead
\subhead 2.1\endsubhead
In this section we fix $w\in W$ and $s\ge1$. Let $\bold1\in\Irr(W)$ be the
unit representation. It is known that if $g\in G$ is unipotent then
$H^i(\cb_g)_{\bold1}$ is $0$ if $i\ne0$ and is $\bbq$ if $i=0$. Using this and the
definition 0.6(a), we see that $\sha(((B^+wB^+)^1)^{F^s})=\z_{w,p^s,\bold1}$. Hence

(a) $\sha((\fB_w^{(1)})^{F^s})=\sha(\cb^{F^s})\z_{w,p^s,\bold1}$. 
\nl
We now write (in $\ch$):
$(-\qq)^{|w|}T_{w\i}\i=\sum_{z\in W}(-1)^{|z|}R_{z,w,\qq}T_z$ where $R_{z,w,\qq}\in\ZZ[\qq]$.
Note that $R_{z,w,\qq}$ is the same as $R_{z,w}$ in \cite{KL, (2.0.a)}
(we use \cite{KL, 2.1(i)}). The following result could be deduced from \cite{KA} (I thank Minh-Tam Trinh
for this remark); but our proof gives additional information, see 2.6(a).

\proclaim{Proposition 2.2} We have $\z_{w,p^s;\bold1}=p^{\nu s}R_{1,w,p^s}$. Hence
$$\sha((\fB_w^{(1)})^{F^s})=\sha(\cb^{F^s})p^{\nu s}R_{1,w,p^s}.$$
\endproclaim

\subhead 2.3\endsubhead
For $y\in W$ we write, as in \cite{LY, 1.3}:
$$a_{w,y}=\sum_{E\in\Irr(W)}\z_{w,p^s;E}tr(y,E)\in\ZZ.$$
From \cite{LY, 2.7} we have:

\proclaim{Proposition 2.4} For any $y\in W$ we have
$$a_{w,y}=\sum_{\ce\in\Irr(W)}
(-1)^{|y|}(\cf_{\ce,s}:\dim(R^1_{y,s})R^1_{y,s})\tr(T_w,\ce_s)p^{\nu s}
\sha(\cb^{F^s})\i.\tag a$$
\endproclaim
Here $(\cf_{\ce,s}:?)$ is the multiplicity of
$\cf_{\ce,s}$ in the virtual representation $?$ of $G^{F^s}$.

\subhead 2.5\endsubhead
We have
$$\sum_{y\in W}a_{w,y}=\sha(W)\z_{w,p^s;\bold1}\tag a$$
hence
$$\align&\z_{w,p^s;\bold1}=\sha(W)\i\sum_{\ce\in\Irr(W)}
\\&(\cf_{\ce,s}:\sum_{y\in W}\dim(R^1_{y,s})D(R^1_{y,s}))
\tr(T_w,\ce_s)p^{\nu s}\sha(\cb^{F^s})\i,\endalign$$
where $D$ is the duality homomorphism from the group of virtual representations of $G^{F^s}$
to itself (see \cite{L84, 6.8}); recall that $D(R^1_{y,s})=(-1)^{|y|}R^1_{y,s}$.
We now replace
$(\cf_{\ce,s}:\sum_{y\in W}\dim(R^1_{y,s})D(R^1_{y,s}))$ by
$(D(\cf_{\ce,s}):\sum_{y\in W}\dim(R^1_{y,s})R^1_{y,s})$ and then by
$(\cf_{\ce\ot\sgn,s}:\sum_{y\in W}\dim(R^1_{y,s})R^1_{y,s})$
where $\sgn$ is the sign representation of $W$. By \cite{DL, 7.5}, 
$\sha(W)\i\sum_{y\in W}\dim(R^1_{y,s})R^1_{y,s})$ is the projection of the regular
representation $Reg$ of $G^{F^s}$ onto the space spanned by unipotent representations.
It follows that
$$\align&\z_{w,p^s;\bold1}=\sum_{\ce\text{ in }\Irr(W)}(\cf_{\ce\ot\sgn,s}:Reg)
\tr(T_w,\ce_s)p^{\nu s}\sha(\cb^{F^s})\i\\&=
\sum_{\ce\text{ in }\Irr(W)}\dim(\cf_{\ce\ot\sgn,s})\tr(T_w,\ce_s)p^{\nu s}\sha(\cb^{F^s})\i.
\endalign$$

It is known that $\tr(T_w,\ce_s)=\tr((-p^s)^{|w|}T_{w\i}\i,(\ce\ot\sgn)_s)$. It follows that
$$\align&\z_{w,p^s;\bold1}=
\sum_{\ce\text{ in }\Irr(W)}\dim(\cf_{\ce\ot\sgn,s})\tr((-p^s)^{|w|}T_{w\i}\i,(\ce\ot\sgn)_s)
p^{\nu s}\sha(\cb^{F^s})\i\\&=
\sum_{\ce'\text{ in }\Irr(W)}\dim(\cf_{\ce',s})\tr((-p^s)^{|w|}T_{w\i}\i,\ce'_s)
p^{\nu s}\sha(\cb^{F^s})\i\\&=
\tr((-p^s)^{|w|}T_{w\i}\i,\cf_s)p^{\nu s}\sha(\cb^{F^s})\i\\&
=\sum_{z\in W}(-1)^{|z|}R_{z,w,p^s}\tr(T_z;\cf_s)p^{\nu s}\sha(\cb^{F^s})\i.\endalign$$  

We have $\tr(T_z,\cf_s)=0$ if $z\ne 1$ and $\tr(T_1,\cf_s)=\sha(\cb^{F^s})$ hence
$\z_{w,p^s;\bold1}=p^{\nu s}R_{1,w,p^s}$. This completes the proof of 2.2.

\subhead 2.6\endsubhead
For any $y\in W$ we set
$$\t_y(w,p^s)=\sum_{E\in\Irr(W)}\tr(T_w,E_s)(\cf_{E,s}:R^1_{y,s})\in\ZZ.$$
Note that $\t_1(w,p^s)=\t(w,p^s)$. (We use that

$\t_1(w,p^s)=\sum_{E\in\Irr(W)}\tr(T_w,E_s)\dim(E)=\t(w,p^s)$.) 
\nl
We have
$$R_{1,w,p^s}(p^s-1)^{\dim\r}=\sha(W)\i\sum_{y\in W}\t_y(w,p^s)\det(p^s-y,\r)\i.
\tag a$$

Using 2.2, 2.5(a),2.4 we have
$$\align&R_{1,w,p^s}=\sha(W)\i\sum_{y\in W}\sum_{\ce\in\Irr(W)}
(-1)^{|y|}(\cf_{\ce,s}:\dim(R^1_{y,s})R^1_{y,s})\tr(T_w,\ce_s)\sha(\cb^{F^s})\i\\&
=\sha(W)\i\sum_{y\in W}\sum_{\ce\in\Irr(W)}
(-1)^{|y|}\dim(R^1_{y,s})\t_y(w,p^s)\sha(\cb^{F^s})\i.\endalign$$
It remains to use the known formula
$$\dim(R^1_{y,s})=(-1)^{|y|}\sha(\cb^{F^s})(p^s-1)^{\dim\r}\det(p^s-y,\r)\i.$$

\head 3. Smoothness\endhead
\subhead 3.1\endsubhead

In this section we fix $w\in W$ which is elliptic and of minimal length in its conjugacy
class. 

\proclaim{Proposition 3.2} (a) $\tfB_w$ is smooth, irreducible, of dimension
$r+2\nu+|w|$.

(b) Let $\d\in\ct$. Then $\tfB_w^\d$ is smooth, irreducible, of dimension $2\nu+|w|$.
\endproclaim

\subhead 3.3\endsubhead
For a point $x$ of a smooth variety $V$ let $T_x(V)$ be the tangent space of $V$ at $x$.

Let $E$ be $\{(g,B)\in G\T\cb; g\in B\}$ (in case 3.2(a)) or
$\{(g,B)\in G\T\cb;g\in B_\d\}$ (in case 3.2(b)).
Then $E$ is a smooth variety of dimension $e$ where $e=r+2\nu$ (in case 3.2(a)) or $e=2\nu$
(in case 3.2(b)).
Let $\tfB_w(E)=\{(g,B,B')\in\tfB_w;(g,B)\in E\}$. This equals $\tfB_w$ (in case 3.2(a)) and
$\tfB_w^\d$ (in case 3.2(b)). We can identify $\tfB_w(E)$ with the intersection

(a) ${}'\tfB_w(E)\cap{}''\tfB(E)$
\nl
where

${}'\tfB_w(E)=\{(g,B,B',B'')\in G\T\cb\T\cb\T\cb;(g,B)\in E,(B',B'')\in\co_w\}$,

${}''\tfB(E)=\{(g,B,B',B'')\in G\T\cb\T\cb\T\cb;(g,B)\in E,gB'g\i=B''\}$,
\nl
are smooth subvarieties of pure dimension $e+\nu+l(w), e+\nu$ of the smooth variety
$\tfB(E)=\{(g,B,B',B'')\in G\T\cb\T\cb\T\cb;(g,B)\in E\}$ of pure dimension $e+2\nu$.
We show that

(b) the intersection (a) is transversal; hence $\tfB_w(E)$ is smooth of pure dimension
$e+l(w)$.    
\nl
Let $x=(g,B,B',B'')\in\tfB_w(E)$ (viewed as a subset of $\tfB(E)$).
Let $V_g=\{(B',B'')\in\cb\T\cb;gB'g\i=B''\}$. From \cite{L11, 5.6} we have

(c) $T_{(B',B'')}(V_g)+T_{(B',B'')}(\co_w)=T_{(B',B'')}(\cb\T\cb)$.
\nl
We have

$T_x(\tfB(E))=T_{g,B}(E)\T T_{(B',B'')}(\cb\T\cb)$,

$T_x({}'\tfB_w(E))=T_{g,B}(E)\T T_{(B',B'')}(\co_w)$,

$T_x({}''\tfB(E))=T_{g,B}(E)\T T_{(B',B'')}(V_g)$,
\nl
hence using (c) we have
$$\align&T_x({}'\tfB_w(E))+T_x({}''\tfB(E))=T_{g,B}(E)\T(T_{(B',B'')}(V_g)+T_{(B',B'')}(\co_w))\\&=
T_{g,B}(E)\T T_{(B',B'')}(\cb\T\cb)=T_x(\tfB(E)).\endalign$$
This proves (b).
\nl
The argument above shows also that

(d) {\it the map $\p:\tfB_w@>>>\ct$ in 0.5 is smooth.} 
\nl
(I thank Zhiwei Yun for this remark.)

\subhead 3.4\endsubhead
It is known that for $w$ elliptic, $\t(w,\qq)$ is a monic polynomial in $\qq$ of degree $|w|$. (This
was stated without proof in \cite{L78, p.27, line 2-4}; a proof is given in \cite{L20}.)
Using this and 0.4(b) we see that for $s=1,2,\do$,
the number of $\FF_{p^s}$-rational points of 
$\tfB_w$ is a monic polynomial in $p^s$ of degree equal to $r+2\nu+|w|$.
Using instead 0.5(a) for $s$ such that $\d$ in 3.2(b) is defined over $\FF_{p^s}$
we see that for such $s$, the number of $\FF_{p^s}$-rational points of 
$\tfB_w^\d$ is a monic polynomial in $p^s$ of degree equal to $2\nu+|w|$.
From this the irreducibility statements in 3.2 follow (we use 3.3(b)). This completes the
proof of 3.2.

\head 4. Quotient by the $G$-action\endhead
\subhead 4.1\endsubhead
In this section we assume that $G$ is semisimple
and that $w\in W$ is as in 0.8. Let $s\ge1$.
Now $G$ acts on $\tfB_w$ by $x:(g,B,B')\m(xgx\i,xBx\i,xB'x\i)$ and
on $\fB_w$ by $x:(g,B')\m(xgx\i,xB'x\i)$.
Let $G\bsl\tfB_w$, $G\bsl\fB_w$ be the set of $G$-orbits on $\tfB_w$, $\fB_w$.   
Now $F^s:\tfB_w@>>>\tfB_w$, $F^s:\fB_w@>>>\fB_w$ induce bijections 
$F^s:G\bsl\tfB_w@>>>G\bsl\tfB_w$, $F^s:G\bsl\fB_w@>>>G\bsl\fB_w$. We show:

(a) $\sha((G\bsl\tfB_w)^{F^s})=\sha((\tfB_w)^{F^s})\sha(G^{F^s})\i$,

(b) $\sha((G\bsl\fB_w)^{F^s})=\sha((\fB_w)^{F^s})\sha(G^{F^s})\i$.

It is enough to show that

(c) if $E$ is any $G$-orbit on $\tfB_w$ such that $F^s(E)=E$
then $\sha(E^{F^s})=\sha(G^{F^s})$;

(d) if $E'$ is any $G$-orbit on $\fB_w$ such that $F^s(E')=E'$
then $\sha(E'{}^{F^s})=\sha(G^{F^s})$.

Now from \cite{L11, 5.2} we see that the isotropy group in $G$ of any point in $E'$ is finite
and hence the isotropy group in $G$ of any point in $E$ is finite. By Lang's theorem $E$ and $E'$
are of the form $G/\G$ with $\G$ a finite subgroup of $G$ stable under $F^s$. We must show that
$\sha((G/\G)^{F^s})=\sha(G^{F^s})$ or that
$\sum_i(-1)^i\tr(F^s,H^i_c(G/\G))=\sum_i(-1)^i\tr(F^s,H^i_c(G))$.
It is enough to show that $H^i_c(G/\G)=H^i_c(G)$ for each $i$. This follows from the fact that
$\G$ acts trivially on $H^i_c(G)$ (the $\G$-action on $G$ by right translations is part of an action
of a connected group).

\subhead 4.2\endsubhead
We show that

(a) $\sha((G\bsl\fB_w)^{F^s})=p^{s|w|}$.
\nl
By \cite{HL} there exists an $F$-stable subvariety $V$ of $G$ stable under
conjugation by elements in a finite abelian $F$-stable subgroup $\D$ of $G$ such that  
any $G$-orbit in $\fB_w$ meets $V$ in exactly one $\D$-orbit and such that $V$ is isomorphic to
$\kk^{|w|}$. It follows that $\sha((G\bsl\fB_w)^{F^s}))$ is equal to the number of
$F^s$-stable $\D$-orbits in $V$ hence is equal to $\sum_i(-1)^i\tr(F^s,H^i_c(\D\bsl V))$.
But $H^i_c(\D\bsl V)$ is the $\D$-invariant part of $H^i_c(V)$ which is $0$ if $i\ne2|w|$ and is
one dimensional with trivial action of $\D$ and with $F^s$ acting as $p^{s|w|}$ if $i=2|w|$.
This proves (a).

\subhead 4.3\endsubhead
From 0.4(b) and 4.1(a) we see that

(a) $\sha((G\bsl\tfB_w)^{F^s})=\t(w,p^s)$.

\subhead 4.4\endsubhead
Let $\d\in\ct^{F^s}$.
Let $G\bsl\tfB_w^\d$ be the set of $G$-orbits on $\tfB_w^\d$.
This is a subset of $G\bsl\tfB_w$ stable under $F^s:G\bsl\tfB_w@>>>G\bsl\tfB_w$.
Using 4.1(c) we see that

(a) $\sha((G\bsl\tfB_w^\d)^{F^s})=\sha((\tfB_w^\d)^{F^s})\sha(G^{F^s})\i$.
\nl
Combining this with 0.5(a) we see that

(b) $\sha((G\bsl\tfB_w^\d)^{F^s})=\t(w,p^s)(p^s-1)^{-r}$

\subhead 4.5\endsubhead
Let $(\d)\in(W\bsl\ct)^{F^s}$.
Let $G\bsl\fB_w^{(\d)}$ be the set of $G$-orbits on $\fB_w^{(\d)}$.
This is a subset of $G\bsl\fB_w$ stable under $F^s:G\bsl\fB_w@>>>G\bsl\fB_w$.
Using 4.1(d) we see that

(a) $\sha((G\bsl\fB_w^{(\d)})^{F^s})=\sha((\fB_w^{(\d)})^{F^s})\sha(G^{F^s})\i$.
\nl
Combining (a) (with $(\d)=(1)$) and 2.2 we see that

(b) $\sha((G\bsl\fB_w^{(1)})^{F^s})=R_{1,w,p^s}(p^s-1)^{-r}$.

\subhead 4.6\endsubhead
Let $H$ be the a regular semisimple conjugacy class in $G$ such that $F^s(H)=H$.
There is a unique conjugacy class $C$ in $W$ such that 
for any $g\in H^{F^s}$ the following holds: for some/any
$B\in\cb$ containing $g$ we have $(B,F^s(B))\in\co_y$ for some $y\in C$.
Let $\d\in\ct$ be the image of $g\in H$ in $B/U_B$ where $B\in\cb$ contains $g$.
The $W$-orbit $(\d)$ of $\d$ is independent of the choice of $g,B$ and is $F^s$-stable.
We have $\sha(H^{F^s})=\sha(G^{F^s})\det(p^s-y,\r)\i$ where $y\in C$ and $\r$ is as in 2.3.
Let $g\in H^{F^s}$. Replacing
$\tr(g,\cf_{E,s})$ by $(\cf_{E,s}:R^1_{y,s})$ in 1.2(a) (with $y$ as above) we obtain
$$\sha(Y_{g,w}^{F^s})=\t_y(w,p^s)$$
(notation of 2.6). We have 
$$\sha((\fB_w^{(\d)})^{F^s})=\sum_{g\in H^{F^s}}\sha(Y_{g,w}^{F^s})=
\sha(H^{F^s})\t_y(w,p^s)=\sha(G^{F^s})\det(p^s-y,\r)\i\t_y(w,p^s).\tag a$$
Using 4.5(a) we deduce
$$\sha((G\bsl\fB_w^{(\d)})^{F^s})=\det(p^s-y,\r)\i\t_y(w,p^s).\tag b$$
In particular if $H$ is such that $C=\{1\}$ then
$$\sha((G\bsl\fB_w^{(\d)})^{F^s})=(p^s-1)^{-r}\t(w,p^s).\tag c$$

\head 5. Symmetry\endhead
\subhead 5.1\endsubhead
In this section we assume that $G$ is semisimple.
Let $v$ be an indeterminate. Let $\tch=\QQ(v)\ot_{\QQ(\qq)}\ch$ where $\QQ(\qq)$ is viewed as a
subfield of $\QQ(v)$ by $\qq\m v^2$. We view $\ch$ as a subset $\tch$ by $h\m1\ot h$.
For any $w\in W$ we set $\ti T_w=v^{-|w|}T_w\in\tch$ and
$\ti\t(w)=v^{-|w|}\t(w,\qq)\in\QQ(v)$, $\ti R_{z,w}=v^{-|w|+|z|}R_{z,w,\qq}\in\QQ(v)$.
Let $\x\m\bar\x$ be the field automorphism of $\QQ(\qq)$ given by $\qq\m\qq\i$. This extends to
a field automorphism $\x\m\bar\x$ of $\QQ(v)$ given by $v\m-v\i$. Let $w\in W$.
We have the following symmetry property.

\proclaim{Proposition 5.2} (a) $\ti\t(w)$ is a polynomial with coefficients in $\NN$ in $(v-v\i)$
hence is invariant under $\bar{}$.

(b) For any $z\in W$, $\ti R_{z,w}$ is a polynomial with coefficients in $\NN$ in $(v-v\i)$
hence is invariant under $\bar{}$.
\endproclaim
We prove (a). For any $y,z\in W$ we have

$\ti T_y\ti T_z=\sum_{z'\in W}c_{y,z}^{z'}\ti T_{z'}$.
\nl
If $y=1$ we have $c_{y,z}^{z'}=1$ if $z=z'$, $c_{y,z}^{z'}=0$ if $z\ne z'$.
If $y\ne1$ we write $y=sy'$ with $|s|=1,|y|=|y'|+1$;  we have

$c_{y,z}^{z'}=c_{y',z}^{sz'}+(v-v\i)dc_{y',z}^{z'}$
where $d=1$ if $|sz'|<|z'|$, $d=0$ if $|sz'|>|z'|$.
\nl
We now see (using induction on $|y|$) that $c_{y,z}^{z'}$ is
a polynomial with coefficients in $\NN$ in $(v-v\i)$.
Since $\ti\t(w)=\sum_{z\in W}c_{w,z}^z$ we see that (a) holds.

We prove (b). From \cite{KL, (2.0.b), (2.0.c)} we have 

$\ti R_{z,y}=\ti R_{sz,sy}$ if $|sz|<|z|,|sy|<|y|$

$\ti R_{z,y}=(v-v\i)\ti R_{sz,y}+\ti R_{sz,sy}$ if $|sz|>|z|,|sy|<|y|$.
\nl
We now see (using induction on $|y|$) that $\ti R_{z,y}$ is
a polynomial with coefficients in $\NN$ in $(v-v\i)$ so that (b) holds.

\proclaim{Corollary 5.3} (a) For $w\in W$ we have $\ov{\t(w,\qq)}=(-\qq)^{-|w|}\t(w,\qq)$.

(b) For $w,z$ in $W$ we have $\ov{R_{z,w,\qq}}=(-\qq)^{-|w|+|z|}R_{z,w,\qq}$.
\endproclaim
Note that (b) is contained in \cite{KL, 2.1(i)}. 

\subhead 5.4\endsubhead
In this subsection we assume that $w$ is as in 3.1.
We set

$\un{\t(w,\qq)}=\t(w,\qq)(\qq-1)^{-r}$, $\un{R_{1,w,\qq}}=R_{1,w,\qq}(\qq-1)^{-r}$.
\nl
The rational functions $\un{\t(w,\qq)}$, $\un{R_{1,w,\qq}}$ take integer values when
$\qq$ is specialized to $p^s$, $s=1,2,\do$ (see 4.4(b), 4.5(b)).
It follows easily (by an argument similar to that in \cite{LY, 1.9}) that

(a) $\un{\t(w,\qq)}\in\ZZ[\qq]$, $\un{R_{1,w,\qq}}\in\ZZ[\qq]$.
\nl
From 5.3 we deduce the following symmetry property
of the polynomials in $\qq$ which calculate
$\sha((G\bsl\tfB_w^1)^{F^s})$, $\sha((G\bsl\fB_w^{(1)})^{F^s})$:

(b) $\ov{\un{\t(w,\qq)}}=\qq^{-|w|+r}\un{\t(w,\qq)}$;
$\ov{\un{R_{z,w,\qq}}}=\qq^{-|w|+r}\un{R_{z,w,\qq}}$.
\nl
(We use that $|w|=r\mod2$.)

\enddocument